\documentclass[oneside,english]{amsart}
\usepackage[T1]{fontenc}
\usepackage[latin9]{inputenc}
\usepackage{amsthm}
\usepackage{amssymb}

\makeatletter
\theoremstyle{plain}
\newtheorem{thm}{Theorem}
\theoremstyle{remark}

\theoremstyle{definition}
\newtheorem{example}[thm]{Example}
\theoremstyle{definition}

\theoremstyle{plain}

\makeatother

\usepackage{babel}

\begin{document}
\title{Fomin--Kirillov algebras}
\author{Leandro VENDRAMIN}
\thanks{This work was supported by Conicet and the Alexander von Humboldt
Foundation. The author thanks Bernard Leclerc for conversations about cluster
algebras and Lie theory.}
\address{Philipps-Universit\"at Marburg\\ 
	FB Mathematik und Informatik \\
	Hans-Meerwein-Stra\ss e\\
	35032 Marburg, Germany}
\email{lvendramin@dm.uba.ar}
\begin{abstract}
	This is an extended abstract of the talk given in the Oberwolfach
	miniworkshop ``Nichols algebras and Weyl groupoids'' in October 2012.
\end{abstract}
\maketitle

For an integer $n\geq3$ denote by $\mathcal{E}_{n}$ the algebra (of type
$A_{n-1}$) with generators $x_{(ij)}$, where $1\leq i<j\leq n$, and relations
\begin{align*}
&x_{(ij)}^{2}=0,&\text{ for $1\leq i<j\leq n$,}\\
&x_{(ij)}x_{(jk)}=x_{(jk)}x_{(ik)}+x_{(ik)}x_{(ij)},&\text{ for $1\leq i<j<k\leq n$,}\\
&x_{(jk)}x_{(ij)}=x_{(ik)}x_{(jk)}+x_{(ij)}x_{(ik)},&\text{ for $1\leq i<j<k\leq n$,}\\
&x_{(ij)}x_{(kl)}=x_{(kl)}x_{(ij)},&\text{ for any discinct $i,j,k,l$.}
\end{align*}

The algebras $\mathcal{E}_n$ are graded by $\mathrm{deg}(x_{(ij)})=1$. Of
course, it is natural to ask if $\mathcal{E}_{n}$ is finite-dimensional. It is
known that $\mathcal{E}_{n}$ is finite-dimensional if $n\leq5$. It has been
conjectured that $\dim\mathcal{E}_n=\infty$ for $n\geq6$.

\begin{example}
The algebra $\mathcal{E}_3$ has dimension $12$. The Hilbert series
$\mathcal{H}_3(t)$ of $\mathcal{E}_3$ is a polynomial of degree $4$:
$\mathcal{H}_3(t)=(2)_{t}^{2}(3)_{t}$, where $(k)_t=1+t+\cdots+t^{k-1}$.
\end{example}

\begin{example}
Computer calculations yield $\dim\mathcal{E}_4=576$. The Hilbert series
$\mathcal{H}_4(t)$ of $\mathcal{E}_4$ is a polynomial of degree $12$:
$\mathcal{H}_4(t)=(2)_{t}^{2}(3)_{t}^2(4)^2_t$.
\end{example}

\begin{example}
Computer calculations yield $\dim\mathcal{E}_5=8294400$. The Hilbert series
$\mathcal{H}_5(t)$ of $\mathcal{E}_5$ is a polynomial of degree $40$:
$\mathcal{H}_5(t)=(4)_{t}^{4}(5)_{t}^{2}(6)_{t}^{4}$.
\end{example}

\begin{example}
The Hilbert series $\mathcal{H}_6(t)$ of $\mathcal{E}_6$ cannot be written as a product of
$t$-numbers. Further, 
\[
\mathcal{H}_6(t)=1+15t+125t^2+765t^3+3831t^4+16605t^5+64432t^6+228855t^7+\cdots
\]
\end{example}
In \cite{MR1667680}, Fomin and Kirillov introduced the algebras $\mathcal{E}_n$
as a new model for the Schubert calculus of a flag manifold.  They proved that
$\mathcal{E}_n$ contains a commutative subalgebra isomorphic the cohomology
ring of the flag manifold.  In \cite{MR2209265}, Bazlov proved that Nichols
algebras provide the correct setting for this model of Schubert calculus. But
what is the relation between the algebras $\mathcal{E}_n$ and Nichols
algebras?  

Let $V_n$ be the vector space with basis $\{v_{(ij)}\mid 1\leq i<j\leq n\}$ and
consider the map $c\in\mathbf{GL}(V_n\otimes V_n)$ defined by
\[
c(v_\sigma\otimes v_\tau)=\chi(\sigma,\tau)v_{\sigma\tau\sigma^{-1}}\otimes v_\sigma,\quad\quad
 \chi(\sigma,\tau)=
 \begin{cases}
 	1 & \text{if $\sigma(i)<\sigma(j)$},\\
  -1 & \text{otherwise,}
\end{cases}
\]
where $\sigma$ and $\tau$ are transpositions, and $\tau=(i\;j)$ with $i<j$.
Since $(V_n,c)$ is a braided vector space, it is possible to consider the
Nichols algebra $\mathfrak{B}(V_n)$.  Bazlov proved that $\mathfrak{B}(V_n)$
contains a commutative subalgebra isomorphic to the cohomology ring of the flag
manifold.  

It is known that $\mathfrak{B}(V_n)=\mathcal{E}_n$ if $3\leq n\leq5$; this was
proved by Milinski and Schneider for $n\leq4$, and by Gra\~na for $n=5$.  It
has been conjectured that $\mathfrak{B}(V_n)$ is quadratic and
$\mathfrak{B}(V_n)=\mathcal{E}_n$; see for example \cite{MR2106930} and
\cite{MR2209265}.

There are many others interesting conjectures about Fomin--Kirillov algebras.
In \cite{MR2106930}, Majid wrote that it might be possible to find a relation
between Fomin--Kirillov algebras and the representation theory of preprojective
algebras. 

Let $\Lambda$ be the preprojective algebra of a quiver of type $A_{n-1}$. It is
known that the number of indecomposable modules over $\Lambda$ is $4$ if $n=3$,
$12$ if $n=4$, and $40$ if $n=5$. Further, $\Lambda$ is of infinite
representation type if $n\geq6$. Majid noticed that the number of
indecomposable modules over $\Lambda$ is equal to the degree of the Hilbert
series of $\mathcal{E}_n$, at least for $3\leq n\leq5$. Majid's conjecture does
not have a precise formulation, but it states that this numerology is not an
accident.

To conclude, we restate Majid's observation in terms of cluster algebras.  Let
$n\geq2$, $G=\mathbf{SL}_n$, and $N$ be the subgroup of upper triangular
matrices with ones in the diagonal.  In \cite{MR2110627}, Berenstein, Fomin and
Zelevinski proved that $\mathbb{C}[N]$, the coordinate ring of $N$, has a
cluster algebra structure. Furthermore, the number of clusters of
$\mathbb{C}[N]$ is given by the following table:
\begin{center}
\begin{tabular}{c|c}
  Lie type of $G$ & Number of clusters\tabularnewline
  \hline
  $A_2$ & $4$ \tabularnewline
  $A_3$ & $12$\tabularnewline
  $A_4$ & $40$\tabularnewline
  others & $\infty$
\end{tabular}
\end{center}

Geiss, Leclerc and Schr\"oer established a relation between the number of
clusters of $\mathbb{C}[N]$ and the number of indecomposable modules over the
preprojective algebra $\Lambda$, see for example \cite{MR2827980}.  This
implies that Majid's observation can be translated into the combinatorial
language of cluster algebras.

\def\cprime{$'$}

\end{document}